\documentclass[a4paper]{amsart}
\usepackage{epsfig}
\usepackage{amsmath}
\usepackage{amsthm}
\usepackage{amssymb}
\usepackage{latexsym,amsfonts}
\usepackage{amscd}
\usepackage[all]{xy}
\parskip\smallskipamount
\numberwithin{equation}{section}
\newcommand{\ck}{{{k}}} 

\newcommand{\mH}{H}

\newcommand{\mand}{\mbox{ and }}

\newcommand{\mor}{\mbox{ or }}

\newcommand{\mfor}{\mbox{ for }}
\newcommand{\mforall}{\mbox{ for all }}

\newcommand{\mHom}{\mbox{Hom}\,}

\newcommand{\Ff}{{\mathbb F}}

\newcommand{\Pp}{{\mathbb P}}

\newcommand{\Rr}{{\mathbb R}}

\newcommand{\cD}{{\cal{D}}}

\newcommand{\cF}{{\cal{F}}}

\newcommand{\cO}{{\cal{O}}}

\newcommand{\cR}{{\cal{R}}}

\newcommand{\cs}{\underline s}

\newcommand{\lra}{\longrightarrow}

\newcommand{\ra}{\rightarrow}

\newtheorem{Prop}{Proposition}[section]
\newtheorem{Thm}[Prop]{Theorem}
\newtheorem{Lemma}[Prop]{Lemma}
\newtheorem{Cor}[Prop]{Corollary}

\let\cal\mathcal

\def\mEnd{\operatorname{End}}
\def\mExt{\operatorname{Ext}}

\def\Hom{\operatorname {Hom}}
\def\Ext{\operatorname {Ext}}

\def\rk{\operatorname {rk}}

\def\id{\operatorname {id}}

\begin{document}
\title{Tilting Bundles on Rational Surfaces and Quasi-Hereditary Algebras}
\author{Lutz Hille, Markus Perling}

\date{\today, This work was supported by DFG priority program 1388 'Representation Theory'}
\maketitle
\begin{abstract}
Let $X$ be any rational surface. We construct a tilting bundle $T$ on
$X$. Moreover, we can choose $T$ in such way that its endomorphism
algebra is quasi-hereditary. In particular, the bounded derived category of
coherent sheaves on $X$ is equivalent to the bounded derived category
of finitely generated modules over a finite dimensional
quasi-hereditary algebra $A$. The construction starts with a full
exceptional sequence of line bundles on $X$ and uses universal
extensions. If $X$ is any smooth projective variety with a full
exceptional sequence of coherent sheaves (or vector bundles, or even
complexes of coherent sheaves) with all 
groups $\mExt^q$ for $q \geq 2$ vanishing, then $X$ also admits a
tilting sheaf (tilting bundle, or tilting complex, respectively)
obtained as a universal   
extension of this exceptional sequence. 
\end{abstract}

\tableofcontents

\section{Introduction}\label{Sintro}

Tilting bundles were first constructed by Beilinson on the projective
$n$--space $\Pp^n$ \cite{Bei}. Later Kapranov obtained tilting bundles on
homogeneous spaces \cite{Kap1}. Moreover, many further examples
are known for 
certain monoidal transformations and projective space bundles
\cite{Orlovsemiorth}. More recently, tilting bundles consisting of
line bundles were investigated by the authors \cite{HP2} and
exceptional sequences on 
stacky toric varieties were constructed by Kawamata \cite{Kaw}. It is
also known that varieties admitting a 
tilting bundle 
satisfy very strict conditions, its Grothendieck group of coherent
sheaves is a finitely
generated free abelian group and the Hodge diamond is concentrated on
the diagonal (in characteristic zero) \cite{BH1}.
However, we are still far from a classification of smooth (projective)
algebraic varieties admitting a tilting bundle. The present note is a
step forward in this direction for algebraic surfaces. The converse of
our main result,
if $X$ is a surface admitting a tilting bundle then it is rational, is
still an open problem.
It can be shown  for many surfaces that tilting bundles cannot exist
using the  
classification (see e.g. \cite{BPV}). 
However, there exist surfaces of general type that have
all the strong properties we need: the canonical divisor has no global
sections, the Grothendieck group of coherent sheaves is finitely
generated and free, and the Hodge diamond is concentrated on the
diagonal.

In this note $X$ is a rational surface over an algebraically closed
field $\ck$. We assume it is smooth and
projective. Some results are valid also for any smooth projective
algebraic variety, however our main interest concerns rational
surfaces. The principal aim is to show, that any rational surface
admits a tilting bundle. The proof is constructive and goes in two
main steps. First, we construct on any rational surface a full
exceptional sequence of line bundles (Section \ref{sectlbds}). This
already 
follows from our previous work \cite{HP2}. Moreover we
show, that in such a sequences there are no $\mExt^2$--groups between
the line bundles. In a second step, we define a universal
(co)extension for such 
sequences, and obtain a tilting bundle. The last step, if we use only
universal extensions, coincide with a
construction known in representation theory of so-called
quasi-hereditary algebras, and is called standardization in
\cite{DR}. 

Since our methods work in a much broader context, we try to be as
general as possible. In fact, the last step, the universal extension
can be defined for any exceptional sequence of complexes (objects in
the corresponding derived category). However, we
obtain a tilting complex only if all higher (that is $\mExt^2$ and
higher groups) do vanish. Otherwise, we obtain at least a partial
result (Theorem \ref{Tgeneralequivalence}). We start with our main
results and then explain 
the strategy of the proof together with the content of this work.
A vector bundle $T$ on an algebraic variety $X$ is called {\sl tilting
bundle} if $\mExt^q(T,T) = 0$ for all $q > 0$ and $T$ generates the
derived category of coherent sheaves on $X$ in the following sense:
the smallest triangulated subcategory of the bounded derived category
$\cD^b(X)$ of coherent sheaves on $X$ containing all direct summands
of $T$ is already $\cD^b(X)$ itself. For further notions of generators
we refer to \cite{BondalVdBergh}.

\begin{Thm}\label{THMmain}
Any smooth, projective rational surface $X$ admits a tilting bundle
$T$ on $X$. 
\end{Thm}

This, in particular, yields an equivalence $\Rr\Hom(T,-)$ between the
bounded derived 
category $\cD^b(X)$ of coherent sheaves on $X$ and the bounded derived category
$\cD^b(A)$ of right modules over the finite dimensional endomorphism
algebra $A$ of $T$.  
In fact we will see in Section \ref{sectconstruction} that we have
many choices to 
construct a tilting bundle $T$. First, we choose a sequence of blow
ups and a standard augmentation (see Definition \ref{defstaugment}) to obtain a full
exceptional sequence of line bundles on $X$. Then we can either use
universal extensions or universal coextensions (we have again a choice
for any $\mExt^1$--block, see the final part in section
\ref{sectunisequences}) to obtain a tilting
bundle. So it is desirable and possible to construct tilting bundles
with further 
good properties. One possibility is to keep the ranks of the
indecomposable direct summands of $T$ small. This needs some detailed
understanding of the non-vanishing extension groups and is based on our
previous work \cite{HP2}. The other way is to
obtain an endomorphism algebra with good homological properties. One
natural choice is to construct a tilting bundle $T$ so that $A =
\mEnd(T)$ becomes  a 
so--called quasi--hereditary algebra. Quasi--hereditary algebras have
very nice and well--understood homological properties (see \cite{DR}
for a short introduction). In particular, there is the subcategory
$\cF(\Delta)$ of
$A$--modules with a $\Delta$--filtration, an additive subcategory closed
under kernels and extensions. Such a choice is also of
interest by a second reason. If we deal with exceptional sequences,
where the higher $\mExt$--groups do not vanish, we obtain a functor
between the derived categories that is not an equivalence. Anyway,
using quasi--hereditary algebras we can at least obtain an equivalence
between certain subcategories (Theorem \ref{Tsurfquasi-her} and
Theorem \ref{Tgeneralequivalence}). 

\begin{Thm}\label{Texcquasi-her}
Let $X$ be a rational surface then there exists a tilting bundle $T$ on
$X$ with quasi-hereditary endomorphism algebra. 
\end{Thm}

Let $\varepsilon$ be any set of objects in an abelian category. Then we
define the subcategory $\cF(\varepsilon)$ as the full subcategory of
objects $M$ admitting a filtration $F^0 = 0 \subseteq F^{1} \subseteq F^{2}
\subseteq, \ldots \subseteq F^{l} = M$ for some integer $l$, so that
for any $0 < i < l$ the quotient $F^{i+1}/F^{i}$ is isomorphic to one object
in $\varepsilon $. We consider this category in the particular case
that $\varepsilon$ is an exceptional sequence. If the abelian category
is the category of finitely generated modules over a finite
dimensional  algebra $A$, an exceptional sequence is also called a set
of standardizable modules (see \cite{DR}). This particular
case we consider in detail in Section \ref{sectquasi-her}. If we
restrict, over a 
quasi-hereditary algebra,  to the set of standard modules $\Delta(1),
\ldots, \Delta(t)$ then these modules form an exceptional sequence and
the category of modules with a $\Delta$--filtration is also called the
{\sl category of good modules} over $A$. This category plays an
important role under the equivalence above. In fact, we can obtain an
equivalence between the subcategory of coherent sheaves
$\cF(\varepsilon)$ admitting a filtration by line bundles in the
exceptional sequence with the well understood category of good modules
over the quasi-hereditary endomorphism algebra. At this point it is
desirable to have small ranks (thus line bundles) for the objects in
$\varepsilon$ to keep the category $\cF(\varepsilon)$ large. In fact
the next result is also constructive, however it needs more background
that we develop only in the last section to formulate it in this
way. Also note, that we can take any full exceptional sequence of line
bundles 
obtained from the Hirzebruch surface by any sequence of standard
augmentations (see Definition \ref{defstaugment}) in the following
theorem. The tilting bundle $T$ is then obtained as a universal
extension of the line bundles in the exceptional sequence.

\begin{Thm}  \label{Tsurfquasi-her}
On any rational surface $X$ there exists a full  exceptional
sequence of line bundles $\varepsilon= (L_{1}, \ldots, L_{t})$ and a
tilting bundle $T$,
so that under the equivalence $\Rr\Hom(T,-)$ between the derived categories above the
category of coherent sheaves $\cF(\varepsilon)$ with a filtration by
the line bundles in $\varepsilon$ is equivalent to the category
$\cF(\Delta)$ of good $A$-modules. Moreover, the functor
$\Rr\Hom(T,-)$ maps $L_i$ to $\Delta(i)$.
\end{Thm}

Now it is desirable to have similar constructions for other
varieties, in particular, in any dimension. We assume $X$ is any
smooth projective (or at least complete) algebraic 
variety and $\varepsilon$ is an exceptional sequence. Then the theorem
above generalizes to this situation. For we need to construct a
quasi-hereditary algebra $A$. It appears as the endomorphism algebra
of the universal extension $\overline E$ of $\varepsilon
$. We discuss this construction in detail in Section
\ref{sectunisequences}. Here we only 
need to know, that there exists such a quasi-hereditary algebra
$A$. The construction is completely parallel to the construction in
\cite{DR} for modules over finite dimensional algebras. Note
that we do not need a full sequence anymore, however the category
$\cF(\varepsilon)$ can be rather small.

\begin{Thm}\label{Tgeneralequivalence}
If $\varepsilon$ is an exceptional sequence of sheaves on $X$, then
there exists a quasi-hereditary algebra $A$ so that the the category
$\cF(\varepsilon)$ of coherent sheaves with an
$\varepsilon$--filtration is equivalent to the category $\cF(\Delta)$ of
$\Delta$--good $A$--modules.    
\end{Thm}

Note that the equivalence in the theorem does, in general, not induce
an equivalence between the corresponding derived categories. To obtain
an equivalence of the derived categories it is necessary (and also
sufficient, as the next theorem shows) that all higher extension
groups vanish. Consequently, we eventually consider exceptional
sequences $\varepsilon$ with vanishing higher extension groups. For
those sequences we can even construct a tilting object. In fact, the
exceptional sequence we start with need not to be a sequence of
sheaves, it can even consist of complexes of sheaves. However, for
complexes we can not expect to get an equivalence for the filtered
objects as above.   For any exceptional sequence
of complexes of coherent sheaves 
$\varepsilon$ we define $\cD(\varepsilon)$ to be the smallest full
triangulated subcategory of $\cD^b(X)$ containing all objects in
$\varepsilon$. In case $\varepsilon$ consists of coherent sheaves the
category $\cF(\epsilon)$ also generates $\cD(\varepsilon)$.

\begin{Thm}\label{ThmUniExt}
Let $\varepsilon = (E_{1}, \ldots, E_{t})$ be any exceptional sequence
of objects in the bounded derived category of coherent sheaves on a
smooth projective algebraic $X$ with $\mExt^{q}(E_{i}, E_{j}) =0$ for
all $q<0$ and all $q \ge 2$. Then the full triangulated subcategory
$\cD(\varepsilon)$ generated by $\varepsilon$ admits a tilting object
$T$, that is 
obtained as a universal extension by objects in $\varepsilon$. If the
exceptional sequence $\varepsilon$  consists of sheaves, then the
tilting object is a sheaf as well, and if the exceptional sequence
consists of vector bundles then $T$ is also a vector bundle. Finally,
if $\varepsilon$ is full then $T$ is a  tilting object in $\cD^{b}(X)$.  
\end{Thm} 

{\em Outline}
The reader just interested in a construction of a tilting bundle on a
rational surface only
needs to read Section 2 to Section 4. The construction gets more
technical if one wants to construct tilting bundles with further
properties or wants to obtain the more general results for any
projective algebraic variety $X$. In fact the results in Section 3 and
Section 4 have a nice interpretation in terms of differential graded
algebras (DG-algebras). Given a DG-algebra as an endomorphism algebra
of an exceptional sequence with only degree zero and degree one terms
(that is $\mExt^q = 0$ for $ \geq 2$) then it is derived equivalent to
an ordinary algebra (that is the endomorphism algebra of the universal
extension). 

In Section \ref{sectlbds} we start with the construction of an
exceptional sequence on any rational surface and prove some further
vanishing results. In section \ref{sectpairs} we consider universal
(co)extensions of pairs of objects. Since we need to use the
construction recursively, it is not sufficient to consider only
exceptionl pairs. In Section \ref{sectunisequences} we define
universal extensions of exceptional sequences and prove Theorem
\ref{THMmain}. In Section \ref{sectquasi-her} we proceed with
quasi-hereditary algebras. In fact we need this notion to define
modules with good filtration. Moreover, using this notion we also get
results for exceptional sequences on varieties of higher
dimension. Finally, in Section \ref{sectconstruction} we construct one
tilting bundle explicitly.

\section{Exceptional sequences of line bundles on rational
  surfaces}\label{sectlbds} 

In this section we construct on any rational surface a full
exceptional sequence of line bundles $\varepsilon = (L_{1}, \ldots,
L_{t})$ that satisfies $\mExt^{2}(L_{i}, L_{j}) = 0$ for all $1 \le
i,j \le t$. The construction is based on a form of a pull back of an
exceptional sequence of line bundles for any blow up (Theorem
\ref{Texistlinesurf}), called standard augmentation in \cite{HP2}.
The main steps of the construction have  
already been proved in \cite{HP2}, Section 5. We only recall the  main
augmentation lemma and prove the vanishing of the higher
Ext-groups. Since any rational surface, not isomorphic to $\Pp^{2}$,
admits a blow down to a Hirzebruch surface in finitely many steps we
can use induction on blow ups. For $\Pp^2$ such a full exceptional
sequence is $(\cO, \cO(1), \cO(2))$. For the Hirzebruch surfaces we
construct an infinite sequence of these sequences (where we can assume
without loss of generality that $L_1$ is the trivial line bundles
$\cO$). Denote by $\Ff_m$ the $m$th Hirzebruch surface. It admits a
natural projection to $\Pp^1$ with fiber $\Pp^1$ and a natural projection
to $\Pp^2$ with exceptional fiber a prime divisor $E$. Then there is a
prime divisor $Q$ linear equivalent to $mP + E$. Computing the
self-intersection numbers we obtain $P^2 = 0, E^2= -m$, and $Q^2 = (mP
+ E)^2 = m$. Moreover, the Picard group of $\Ff_m$ is freely generated
by $\cO(P)$ and $\cO(Q)$ (or $\cO(E)$, respectively), so any line
bundle is isomorphic to 
$\cO(\alpha P + \beta Q)$. The computation of the cohomology groups
uses standard formulas in toric geometry (see \cite{FultonBook},
Section 3.5). The
concrete results for the Hirzebruch surfaces can also be found in
\cite{HGoett}. 

\begin{Prop} \label{THMexcHirzebruch}
Any Hirzebruch surface $\Ff_{m}$ admits a full exceptional
sequence of line bundles $\varepsilon = (\cO, \cO(P), \cO(Q + aP),
\cO(Q + (a+1)P))$, 
where this sequence is strong precisely  when $a \ge 0$.    
\end{Prop}

{\em Proof.}
From the standard formula for cohomology of line bundles on toric
varieties (see \cite{FultonBook} or \cite{OdaBook}) follows $\mfor \beta = -1,$
 or $\beta \geq 0 \mand \alpha \geq -1$  
$$
H^{{1}}(\Ff_{m}, \cO(\alpha P + \beta Q)) = 0 
$$
For line bundles in the the sequence above, the second cohomology
group vanishes, since for any $\beta \geq -1$
$$
H^{2}(\Ff_{m}, \cO(\alpha P + \beta Q)) = 0. 
$$
Consequently, we have on any Hirzebruch surface an infinite family of 
exceptional sequences and an infinite family of strongly
exceptional sequences. It remains to show that both are full. To show
this claim it is sufficient that one exceptional sequence in the
family is full: consider the exact sequence 
$$
0 \lra \cO \lra \cO(P)^2 \lra \cO(2P) \lra 0.
$$
Then we can recursively show that $(\cO, \cO(P), \cO(Q + aP), \cO(Q +
(a+1)P))$ is full, precisely when $(\cO, \cO(P), \cO(Q + (a-1)P),
\cO(Q + aP))$ is full, just tensor the sequence above with $\cO(Q + (a-1)P)$.

It is well-known that the sequence above is full for $\Pp^1 \times
\Pp^1$ and the first Hirzebruch surface $\Ff_1$. Moreover, using the
projection $\Ff_m \lra \Pp^1$ and results in \cite{Orlovsemiorth} we see
that $\varepsilon$ is full on any Hirzebruch surface.
\hfill $\Box$
\medskip

{\sc Remark.} On a Hirzebruch surface with $m\ge 3$ the set of
sequences in Theorem \ref{THMexcHirzebruch} is already 
the complete classification of (strongly) exceptional sequences of
line bundles (up to a twist with a line bundle). For $m=0,1,$ and $2$
there is a finite number of further sequences (see
\cite{HGoett}). 
\medskip

In the next step we consider blow ups $\widetilde X \lra X$ in one
point $x \in X$ with exceptional divisor $E \subset X$. If $L$ is a
line bundle on $X$, 
then we denote the pull back of $L$ under a blow up of $X$ with the
same letter. Since $H^q(X; \widetilde L) = H^q(X; L)$ for any line bundle
$L$ on $X$ this notation is convenient and does not lead to any
confusion if we compute extension groups. Let $X$ be 
any surface with an exceptional sequence $\varepsilon = (L_{1},
\ldots, L_{t})$ of line bundles on $X$. Then we obtain a sequence
$\widetilde \varepsilon := (L_{1}(E), \ldots, L_{i-1}(E), L_{i},
L_{i}(E), L_{i+1}, \dots, L_{t})$ on the blow up $\widetilde X$. 
\medskip

{\sc Definition.} \label{defstaugment}
Given an exceptional sequence $\varepsilon$ on $X$. We call the sequence
$\widetilde \varepsilon^{(i)} := (L_{1}(E), \ldots, L_{i-1}(E), L_{i},
L_{i}(E), L_{i+1}, \dots, L_{t})$ on the blow up a {\sl standard
  augmentation} of $\varepsilon$ (at position $i$).
\medskip

Note that we can choose any $i$ to obtain a standard augmentation, so
for any blow up we have $t$ choices to produce a new sequence.
We will show that
$\widetilde \varepsilon$ is exceptional for each $i=1,\ldots,t$. If
$\varepsilon$ is strongly exceptional then in some cases the new
sequence is even strongly exceptional. In this case, there must exist
a section in $\mHom(L_{j}, L_{k})$  for all $j\le i$ and all $k\ge i$
that does not vanish in $x$ (see \cite{HP2}, Theorem 5.8 and the
proof). However, in general 
the new sequence is only exceptional. The more detailed analysis of
when $\widetilde \varepsilon$ is strongly exceptional is needed only
for the concrete construction of the tilting bundle $T$. So we
leave this part to the end in Section \ref{sectconstruction}. 

\begin{Prop}\label{THMaugmentation}
Let $X$ be a surface with an exceptional sequence $\varepsilon$. Then
$\widetilde \varepsilon$ is an exceptional sequence on the blow up
$\widetilde X$ in one point. If $\varepsilon$ is full then $\widetilde
\varepsilon$ is also full. 
\end{Prop}

{\sc Proof.}
We prove first the vanishing result. For we have to show that on
$\widetilde X$ 
$$
\mExt^{q}(L_{j}, L_{k}(E)) = 0 = \mExt^{q}(L_{j}, L_{k}) \mforall q
\mand j > k.
$$
The second vanishing is obvious, since it coincides with the same
group on $X$. It remains to show the first vanishing.  
Using $\mExt^{q} (L_{j}, L_{k}(E)) = \mH^{q}(L_{j}^{-1}\otimes
L_{k}(E))$ it is sufficient to show the following lemma. Moreover,
also $\mExt^q(L_i(E),L_i) = H^q(\cO(-E)) = 0$ for $q=0,1,2$ follows
from the exact 
sequence in the proof of the following lemma.

Finally, we need to show that $\widetilde \varepsilon$ is full,
provided $\varepsilon$ is. For we consider the semi--orthogonal
decomposition of $\cD^{b}(\widetilde X)$ with respect to $\cD^{b}(X)$
and $\cO_{E}(-1)$ (see \cite{Orlovsemiorth}). Obviously, by
assumption, the line bundles $L_{i}$ 
generate $\cD^{{b}}(X)$. Moreover, the bundles $L_{i}(E)$ and $L_{i}$
generate  $\cO_{E}(-1)$. Then in the last step, the subcategory
generated by $L_{j}(E)$ and  $\cO_{E}(-1)$ contains
$L_{j}$. Consequently, $\widetilde \varepsilon$ generates  $\cD^{b}(X)$ and
$\cO_{E}(-1)$, thus also   $\cD^{b}(\widetilde X)$. 
$\Box$
\medskip

\begin{Lemma}\label{LExt2}
If $L$ is a line bundle on a surface $Y$ with $\mH^{q}(Y; L)= 0$ for all
$q$ and $E \simeq \Pp^{1}$ is a $(-1)$--curve on $Y$ so that
$L\mid_{E}$ is trivial then   $\mH^{q}(Y; L(E))= 0$ for all $q$.   
\end{Lemma}

Note that this is exactly our situation. If we consider the pull back of a
line bundle $L$ to a blow up, then the restriction of $L$ to the
exceptional divisor is trivial and the exceptional divisor is a $(-1)$--curve.

{\sc Proof. }
We consider the short exact sequence
$$
0 \longrightarrow \cO_{\widetilde X} \longrightarrow
\cO_{\widetilde X}(E) \longrightarrow 
\cO_{\Pp^{1}}(-1) \longrightarrow 0 
$$
and tensor it with $L$. Then in the corresponding long exact sequence
the following groups vanish: $\mH^{q}(Y; L)$ for all $q$ and
$\mH^{q}(\Pp^{1}, L\mid_{E}(-1)) =    \mH^{q}(\Pp^{1},
\cO_{\Pp^{1}}(-1)) $ for all $q$. Consequently, the claim follows. 
$\Box$
\medskip

Our main theorem in this section states that any rational surface $X$
admits a full exceptional sequence of line bundles, so that all the
groups  $\mExt^{2}(L_{i}, L_{j})$ vanish. Such a sequence can be
obtained by recursive standard augmentation from an exceptional
sequence on a Hirzebruch surface. 

\begin{Thm}\label{Texistlinesurf}
Let $X$ be any rational surface. Then $X$ admits full exceptional
sequences of line bundles, obtained from a full exceptional sequence
of line bundles on a Hirzbruch surface by applying any standard
augmentation in each step of the blow up. For such a sequence the
groups $\mExt^{2}$ between any two members of the sequence vanish. 
\end{Thm}

We will see in the next two sections that any exceptional sequence with this
property defines a tilting bundle. So, using the universal extension
(to define in the next sections)
we have proved Theorem \ref{THMmain}.

{\sc Remark. }
Note that we defined standard augmentation in \cite{HP2} in a more
general sense and allowed to blow up several times in one step. This
gave us even more flexibility in constructing exceptional
sequences. However, if we perform a standard 
augmentation only for the blow up of one point (as we do in this note)
then it is always
admissible in the sense of  \cite{HP2}, Section 5, and it is
sufficient to prove our results.
\medskip

{\sc Proof.} Since we have already shown the existence of a full
exceptional sequence for any recursive blow up of a Hirzebruch
surface (Proposition \ref{THMaugmentation})  we get a full exceptional
sequence on 
any rational surface $X$. Then, using only recursive standard augmentations, we
obtain a full exceptional sequence $\varepsilon = (L_1,\ldots,L_t)$
with $\mExt^2(L_i,L_j) = 0$ for all $1 \leq i,j \leq t$: just apply
the proof of Lemma \ref{LExt2} also to the exact sequence
$$
0 \longrightarrow \cO_{\widetilde X}(E) \longrightarrow
\cO_{\widetilde X} \longrightarrow 
\cO_{\Pp^{1}} \longrightarrow 0.
$$
Then $\widetilde \varepsilon$ has no $\mExt^2$ between any members of
the sequence since $H^2(\widetilde X; L(-E)) = 0$ and thus also
$\mExt^2(L_i(E),L_j) = 0 = \mExt^2(L_i,L_j)$. 
\hfill $\Box$

For later use in Section \ref{sectconstruction} we also need to
prove some further vanishing results. In particular, we need to
compute the cohomology of $\cO(R_i - R_j)$, where $R_i$ and $R_j$ are
both divisors (not prime in general) of self-intersection $-1$. So let
$X$ be any rational surface, not $\Pp^2$  and fix a sequence of blow downs of a
$(-1)$--curve in $X_i$ step by step to a Hirzebruch surface
$$
X = X_{t} \lra X_{t-1} \lra \ldots \lra X_1 \lra X_0 = \Ff_m.
$$
Note that the rank of the Grothendieck group of $\Ff_m$ is $4$ so for
$X$ it is just $t+4$, where $t+4$ is also the length of the full
exceptional sequence. In each $X_i$ (for $i > 0$) we a distinguished
$(-1)$-curve $E_i$ blown down under $X_i \lra X_{i-1}$. We define
$R_i$ to be the divisor on $X$ obtained as the pull back of $E_i$ to
$X$. Then we need to compute $H^q(X, \cO(R_i - R_j))$.

{\sc Definition. }
We say $R_i$ is above $R_j$ if $E_i$ is blown down to a point $P_i \in
X_{i-1}$ that is on the inverse image of $E_j$ in $X_{i-1}$. Then we also
write $i \succ j$.

If $R_i$ is above $R_j$ then $i$ must be larger than $j$,

\begin{Lemma}\label{LExtOR}
The cohomology groups $H^1(X, \cO(R_i - R_j))$ and  $H^0(X, \cO(R_i - 
R_j))$ both equal $\ck$ precisely
when $i \succ j$, otherwise both vanish. The second cohomology group
$H^2(X, \cO(R_i - R_j))$ is always zero.
\end{Lemma}

{\sc Proof. }
First note that  $H^2(X, \cO(R_i - R_j)) =  H^2(X_i, \cO(E_i - R_j))$
vanishes, since $H^2(X_j, \cO(-R_j))  = H^2(X_i, \cO(-R_j))$ is zero and
we have Proposition \ref{THMaugmentation}. Further note that by
Riemann-Roch the Euler characteristic of $\cO(E_i - R_j)$ is zero,
thus $H^0(X, \cO(R_i - R_j))$ and $H^1(X, \cO(R_i - R_j))$ have the same
dimension. \\
We start to compute  $H^0(X, \cO(R_i - R_j))$. If $R_i$ is above $R_j$
then $R_i - R_j$ is effective, thus  $H^0(X, \cO(R_i - R_j)) \not=
0$. On the other hand  $H^0(X, \cO(R_i)) = H^0(X_i, \cO(E_i))$ is
one-dimensional and the dimension of $H^0(X, \cO(R_i - R_j))$ can not
exceed the dimension of $H^0(X, \cO(R_i)$. A similar argument applies
to $H^0(X, \cO(R_i - R_j))$ for $R_i$ not above $R_j$, then we obtain
$H^0(X, \cO(R_i - R_j)) = 0$. \hfill $\Box$

\section{Universal extension of pairs}\label{sectpairs}

In this section we study universal extensions of objects in an abelian
category. We work with sheaves, however all the techniques developed
here work whenever the groups $\mExt^1$ are finite dimensional. The
principal aim of this part is to show, that any pair of objects
with certain $\mExt$--groups vanishing can be transformed,
using universal (co)extensions, in a pair with vanishing
$\mExt^1$--group. Roughly spoken we apply   a certain partial mutation
(in the sense of \cite{BondalAssocAlg}) to such a pair and the first
extension group vanishes. However, the 
price we have to pay is, that we create new homomorphisms between the
new objects. In particular, whenever we have a non-vanishing extension
group we create, using universal extensions, additional homomorphism
in the opposite direction. Thus the result is no longer an exceptional
sequence. Even worse, the new object has nontrivial endomorphisms.
\medskip


{\sc Definition. }
Let $(E,F)$ be a pair of coherent sheaves. 
Then we define the {\sl universal extension} $\overline E$ of $E$
by $F$ (respectively, the {\sl universal coextension} $\overline F$ of $F$
by $E$) by the
following extension 
$$
0 \longrightarrow F \otimes \mExt^{1}(E,F)^{\ast} \longrightarrow
\overline E \longrightarrow E \longrightarrow 0, 
$$
respectively
$$ 
0 \longrightarrow F \longrightarrow \overline F \longrightarrow E
\otimes \mExt^{1}(E,F))\longrightarrow 0. 
$$
If we consider the second exact sequence as a triangle in the derived
category, the boundary map is the canonical evaluation map
$\mExt^1(E,F)\otimes E 
\lra F[1]$. This map induces, just by taking the adjoint, a canonical
map $ E \lra \mExt^1(E,F)^{\ast}\otimes F[1]$. The mapping cone over
this map defines the first exact sequence. Thus, both exact sequences
are unique up to isomorphism. The first exact sequence can be
characterized by the following property: if we apply
$\mHom(-,F)$ in the  long exact sequence we get an
induced map $\mHom(F,F) \otimes \mExt^1(E,F) \lra
\mExt^1(E,F)$ that is the Yoneda product. This map is surjective,
since $\id \otimes \xi$ maps to $\xi$. In a similar way one can
characterize the second exact sequence. If we apply $\mHom(E,-)$ then
we get a surjective map $\mHom(E,E) \otimes \mExt^1(E,F) \lra
\mExt^1(E,F)$ that is just the ordinary Yoneda product.
The following lemma is formulated in more generality than actually
needed.

\begin{Lemma}\label{Lextensionvanishes}
a) Let $(E,F)$ be a pair of objects with $\mExt^{q}(F,E)= \mExt^q(E,E)
= \mExt^q(F,F) = 0$ for all $q
> 0$.  
Then $(E,\overline F)$ and $(\overline E, F)$ satisfy
$\mExt^{1}(E, \overline F) = 0 = \mExt^{q}(\overline F, E)$ and
$\mExt^{1}(\overline E, F) = 0 = \mExt^{q}(F, \overline E)$ 
 for all $q>0$. Moreover, $\mExt^{1}(\overline E, \overline E) = 0 =
 \mExt^{1}(\overline F, \overline F)$. \\ 
b) If in addition we have $\mExt^{q}(E,F)=0$ for some $q\ge 2$ then
$\mExt^{q}(E,\overline F) = 0 = \mExt^{q}(\overline F, E)$ and
$\mExt^{q}(F,\overline E) = 0 = \mExt^{q}(\overline E, F)$. Moreover,
$\mExt^{q}(\overline E, \overline E) = 0 = 
 \mExt^{q}(\overline F, \overline F)$.\\

\end{Lemma}

Roughly spoken we can replace any exceptional pair with only
non-vanishing $\mHom$ and $\mExt^{1}$ by a pair with only
non-vanishing $\mHom$. If we perform this universal extension
recursively we can replace any full exceptional sequence with
vanishing $\mExt^{q}$ for $q>1$ by a tilting object (for the details see Section \ref{sectunisequences}). 

{\sc Proof}. The proof is just a standard diagram chasing, we only
prove the crucial step for a universal extension. 
The proof is analogous for coextensions.
Note that the vanishing of $\mExt^{q}(\overline E, F)$ follows
from the long exact sequence for $\mHom(-,F)$ applied to the universal
extension of $F$ by $E$. We show that $\mExt^{1}(F, \overline E)$
vanishes. We apply $\mHom(-,F)$ to the universal extension. Since the
boundary map $\mHom(F,F) \otimes \mExt^1(E,F) \lra
\mExt^1(E,F)$ is surjective and $\mExt^1(F,F) = 0$ we obtain
$\mExt^{1}(\overline E, F) = 0$. To obtain $\mExt^{1}(\overline E,
\overline E) = 0 $ we apply $\mExt^{1}(\overline E,-)$ to the
universal extension. Since $\mExt^{1}(\overline E,E) = 0$ by
assumption and $\mExt^{1}(\overline E,F) = 0$ by the previous argument
we obtain the desired vanishing.
This finishes the proof of a) for the
universal extension. \\
To show b) we only need to apply $\mHom(-,F)$ to the universal
extension and get the vanishing form the long exact sequence. For the
second vanishing we apply $\mExt^{q}(\overline E,-)$ to the
universal extension, as we did for $q =1$ in a).
\hfill $\Box$
\medskip

{\sc Remark.} Note that $\overline E$, respectively $\overline F$ need
not to be indecomposable. However, the above vanishing result holds
for any direct summand as well. If we deal with an exceptional pair $(E,F)$,
then there exists a unique indecomposable direct summand $\overline
E_1$ of  $\overline E$ with $\dim \mHom(\overline E_1,E) = \dim
\mHom(\overline E,E) = \dim \mHom(E,E) = 1$. This allows later to
distinguish certain indecomposable direct summands of our universal
(co)extension. 

The following lemma is useful for a construction of a tilting bundle
with small rank. In fact we will show in Section
\ref{sectconstruction} that we can on 
any rational surface construct full exceptional sequences of line bundles
with $\dim \mExt^1(L_i,L_j) \leq 1$ for all $i,j$ (there are at most
one-dimensional extension groups).

\begin{Lemma}
Assume $(E,F)$ is a pair of coherent sheaves on $X$ with $\dim
\mExt^1(E,F) = 1$, then the universal extension $\overline E$ is
isomorphic to the universal coextension $\overline F$.
\end{Lemma}

{\sc Proof. }
In this case the end terms of the two defining short exact sequences
coincide. Any non-trivial element $\xi \in \mExt^1(E,F) = \ck$ defines
the same middle term. Since $\overline E$ and $\overline F$ are both
unique up to isomorphism, they must be isomorphic.
\hfill $\Box$
\medskip


\section{Universal extensions of exceptional sequences}
\label{sectunisequences}  

In this section we start with an exceptional sequence and perform
universal extension or coextensions recursively. We explain the
construction for universal extensions, for universal coextensions the
construction is dual. 
\medskip 

{\sc Definition. }
Let $\varepsilon =
(E_1,\ldots, E_t)$ be any exceptional sequence, then we define
$E_i(1) := E_i(2) := \ldots := E_i(i) := E_i$ and $E_i(j)$ to be the
universal extension of 
$E_{i}(j-1)$ by $E_{j}$ for $j > i$. Thus we have exact sequences
$$
0 \lra E_j \otimes \mExt^1(E_i(j-1),E_j)^{\ast} \lra E_i(j) \lra
E_i(j-1) \lra 0 
$$
for all $t \geq j > i$ (for $i \leq j$ the sequences
are always trivial, since the first term vanishes). Thus we have defined
new objects $E_i(t)$ for $1 \leq i \leq t$ 
that are not necessarily indecomposable. We define $\overline E_i$ to
be the unique indecomposable direct summand of $E_i(t)$ with the
property that $\ck = \mHom(E_i(t), E_i) = \mHom(\overline E_i, E_i)$
and denote by $\overline E$ the direct sum $\oplus_{i=1}^t \overline
E_i$. We call $\overline E$ the {\em universal extension} of the exceptional
sequence $\varepsilon$ and $\overline E_i$ the universal extension of
$E_i$ by $E_{i+1},\ldots,E_t$. 

In a dual way we define universal coextensions $\underline E$ and
$\underline E_i$ (here we need to perform recursive universal
coextensions with $E_{i-1}, \ldots, E_1$ instead).

{\sc Remark. }
Note that we can also split the exceptional sequence into two subsets
and perform universal extensions in the first and universal
coextensions in the second subset. 

\begin{Thm}
Let $\varepsilon = (E_1,\ldots, E_t)$ be an exceptional sequence, $E:=
\oplus E_i$ the direct sum of the elements in $\varepsilon$ and
$\overline E$ the direct sum of the objects $\overline E_i$
constructed above. Then $\mExt^1(\overline E, \overline E) = 0$. If,
moreover, $\mExt^q(E,E) = 0$ for some $q > 1$ then also
$\mExt^q(\overline E, \overline E) = 0$. 
\end{Thm}

{\sc Proof. }
We first show $\mExt^1(E_i(j), E_l) = 0$ for all $j \geq l$. Assume
$j = l$ then this follows from Lemma \ref{Lextensionvanishes}, since
$$
0 \ra E_j \otimes \mExt^1(E_i(j-1),E_j)^{\ast} \ra E_i(j) \ra
E_i(j-1) \lra 0 
$$
is a universal extension by $E_j$. Then use induction over $j$, the
induction step just follows from applying $\mExt^1(-,E_l)$ to the
defining exact sequence above.\\
As a consequence we obtain $\mExt^1(E_i(t), E_j) = 0$ for all $1 \leq
j \leq t$. Now we apply $\mExt^1(E_i(t),-)$ to the defining sequence
for $E_i(j)$ and obtain an exact sequence
$$
\mExt^1(E_i(t), E_j) \otimes \mExt^1(E_i(j-1),E_j)^{\ast} \lra
\mExt^1(E_i(t), E_i(j)) \lra 
\mExt^1(E_i(t), E_i(j-1)).
$$
The first and the last term are zero by induction, thus the middle
term vanishes for all $j\geq i$, in particular, it vanishes for $j=t$.
In a similar way we can, using $\mExt^q(\overline E,-)$, show the last
claim.
\hfill $\Box$
\medskip

\begin{Thm}
Assume $\varepsilon = (E_1,\ldots, E_t)$  is a full exceptional
sequence of sheaves (or complexes of sheaves), $E = \oplus_{i=1}^t
E_i$, with $\mExt^q(E,E) = 0$ for all 
$q \not= 0,1$. Then the universal extension $\overline 
E$ is a tilting sheaf (or tilting object).
\end{Thm}

{\sc Proof. }
Using the defining exact sequences we see that $E$ and $\overline E$
generate the same subcategory in the derived category $\cD^b(X)$ (or
the corresponding triangulated category of the abelian category we
work with). Moreover, $\mExt^q(\overline E, \overline E)= 0$ for all
$q \not= 0$ by the previous result. 
\hfill $\Box$
\medskip

{\cs Proof. (of Theorem \ref{ThmUniExt})}
The tilting object is just the universal extension of the exceptional
sequence $\varepsilon$.  
\hfill $\Box$
\medskip

Now we are interested in sequences where we can also perform both,
universal extensions and universal coextensions. To explain this we
collect our objects into $\mExt^1$--blocks. Let $\varepsilon =
(E_1,\ldots, E_t)$ be an exceptional sequence. 
\medskip

{\sc Definition. }
We define a graph
$\Gamma(\varepsilon)$, called {\sl Ext-graph} 
with vertices $i$ for $1 \leq i \leq t$ (the indices of the objects in
$\varepsilon$) and an edge between $i$ and $j$, whenever
$\mExt^1(E_i,E_j)$ or $\mExt^1(E_j,E_i)$ does not vanish. Then an
{\sl $\mExt$--block} consists of a connected component in
$\Gamma(\varepsilon)$. 
\medskip

In this way we define certain subsets, for each
subset we can choose either a universal extension or a universal
coextension. Note that we get different results only if at least one
subset consists of at least three elements or there are
two-dimensional extension groups. 

In the last section we discuss in more detail how to choose an
exceptional sequence of line bundles on a rational surface $X$
depending on the sequence of blow ups from a Hirzebruch surface (or
the projective plane). Then we also know precisely the
$\mExt^1$--blocks. We can use this to minimize the non-vanishing
extension groups so that $\overline E$ has a small rank.


\section{Quasi-hereditary algebras} \label{sectquasi-her}

Let $X$ be a rational surface, or even any algebraic variety, with a tilting
sheaf that is the universal extension of a full exceptional sequence
(with all higher Ext--groups vanishing)
on $X$. Then we claim that the endomorphism algebra of this sequence
satisfies a well-understood and extensively studied property: it is
quasi-hereditary. Note that a quasi-hereditary algebra is an algebra
with an order on its primitive orthogonal idempotents (or equivalently
on its isomorphism classes of indecomposable projective modules). This
order is just induced from the natural order in the exceptional
sequence we started with. 

If we used universal coextensions then we get the dual notion of
so-called $\nabla$--modules. If we use both, universal extensions in
some $\mExt^1$--blocks and universal coextensions in the remaining
$\mExt^1$--blocks, then the endomorphism algebra is not
quasi-hereditary (except the blocks only consist of two members). 
\medskip

In this section we review some of the main properties on
quasi-hereditary algebras. In particular, we use the so-called
standardization introduced by Dlab and
Ringel in \cite{DR} to prove Theorem \ref{Texcquasi-her} and, more
important, Theorem \ref{Tgeneralequivalence}.

Let $A$ be the endomorphism algebra of a sheaf $T$ that is obtained as
a universal extension of an exceptional sequence $\varepsilon =
(E_1,\ldots,E_t)$ of sheaves on
$X$. We decompose $T$ into indecomposable direct summands $T =
\oplus_{i=1}^t T(i)$. Then the natural order in $\varepsilon$ defines
an order on the 
indecomposable projective $A$-modules $P(i) :=
\mHom(T,T(i))$. Moreover, we define $\Delta(i)$ to be the quotient of
$P(i)$ by the maximal submodule generated by any direct sum $\oplus_{j <
  i}P(j)^{a(j)}$. This submodule is a proper submodule.

{\sc Definition. } The algebra $A$ is called
{\sl quasi-hereditary} (with this order) if each $P(i)$ is in
$\cF(\Delta)$ (see \cite{DR}). 

{\sc Remark. }
Note that our definition is slightly different to the one in
\cite{DR}, since an exceptional sequence of sheaves on a complete
variety $X$ has all the properties
of a standardizable set. 

\begin{Thm}
Let $\varepsilon = (E_1,\ldots,E_t)$ be any exceptional sequence of
sheaves and $T$ the 
sheaf obtained from $\varepsilon$ by its universal
extension. Then the endomorphism algebra of $T$ is quasi-hereditary
with $\Delta$--modules $\Delta(i) = \Hom(T,E_i)$.
\end{Thm}

{\sc Proof.}
We need to show that any finitely generated projective $A$-module has
a filtration by the 
modules $\Delta(i)$ defined as a quotient of   $P(i)$. Note that the
recursive universal extensions provides us with such a filtration for
the objects $\overline E_i$ by induction: we consider the defining
exact sequence (from the previous section)
$$
0 \lra E_j \otimes \mExt^1(E_i(j-1),E_j)^{\ast} \lra E_i(j) \lra
E_i(j-1) \lra 0. 
$$
If we apply $\mHom(T,-)$ we obtain an exact sequence of $A$-modules
$$
0 \ra \mHom(T,E_j) \otimes \mExt^1(E_i(j-1),E_j)^{\ast} \ra
\mHom(T,E_i(j)) \ra 
\mHom(T,E_i(j-1)) \ra 0. 
$$
(it is exact, since $\mExt^1(T,E_j) = 0$). Thus, by induction over $j$
the right $A$--module $\mHom(T,E_i(t))$ admits a filtration by the modules
$\Delta(i)$. Now we use that $\cF(\Delta)$ is closed under direct
summands (see e.g.~the characterization in \cite{DR}, Theorem
1). 
\hfill $\Box$
\medskip

\begin{Thm}\label{Thmabove}
Let $\varepsilon = (E_1, \ldots,E_t)$ be a full exceptional sequence
of sheaves and $T$ its universal extension. Then the functor
$\mHom(T,-)$ induces 
an equivalence between $\cF(\varepsilon)$ and $\cF(\Delta)$ mapping
$E_i$ to $\Delta(i)$.
\end{Thm}

{\sc Proof. } 
Let $F$ be any sheaf in  $\cF(\varepsilon)$. Then $\Ext^1(T,F) = 0$
since $\mExt^1(T,E_i) = 0$ for all $i$. Using the exact sequences
defining the filtration of $F$ recursively, we get a filtration of
$\mHom(T,F)$ by $\Delta(i) = \mHom(T,E_i)$. Thus $\mHom(T,F)$ is in
$\cF(\Delta)$. Conversely, let $M$ be an $A$-module in $\cF(\Delta)$,
then $M$ has a projective presentation $P^1 \stackrel{f}{\lra} P^0 \lra M \lra
0$. This defines an induced map $T^1 \stackrel{f^+}{\lra} T^0$, where
$T^1$ and $T^0$ are direct sums of direct summands of $T$. The map
$f^+$ is just defined using the equality $A = \mHom_A(A,A) = \mHom_X(T,T)$ and
the fact that $P^1$ and $P^0$ are direct sums of direct summands of
$A$. Define $F(M)$ to be the cokernel of $f^+$.  Note that $f$ is
injective precisely when $f^+$ is injective. Thus the
$\Delta$-filtration of $M$ induces an filtration of $F(M)$ showing
$F(M)$ is in $\cF(\varepsilon)$. Note that under the functor $M
\mapsto F(M)$ the module $\Delta(i)$ is mapped to $E_i$. Thus this functor is
inverse to $\mHom(T,-)$ finishing the proof. 
\hfill $\Box$
\medskip

{\sc Proof.} (of Theorem \ref{Tgeneralequivalence}) \\
The proof of Theorem \ref{Thmabove} provides us with an explicit
construction for the algebra $A$ as the endomorphism algebra of a
universal extension $T$ of an exceptional sequence.
\hfill $\Box$
\medskip

{\sc Proof.} (of Theorem \ref{Texcquasi-her})\\
If $\varepsilon = (E_1,\ldots,E_t)$ is any exceptional sequence with
$\mExt^2(E_i,E_j) = 
0$ for $1 \leq i,j \leq t$ on a surface $X$, then its universal
extension has a quasi-hereditary endomorphism algebra by the theorem
above. Such a sequence, consisting even of line bundles, exists by
Theorem \ref{Texistlinesurf}.
\hfill $\Box$
\medskip

{\sc Proof.} (of Theorem \ref{Tsurfquasi-her}) \\
This result was proved above, where we replace any sheaf $E_i$ just by
a line bundle $L_i$.
\hfill $\Box$
\medskip

{\sc Remark. }
The principal idea of the theorem above can be found in
\cite{DR}, 3. standardization. Therein is a similar construction for
any abelian category. In fact, such a
construction, even in greater generality, can be performed in any
abelian $\ck$-category with finite dimensional extension groups. 

\section{Construction of tilting bundle on rational
  surfaces}\label{sectconstruction}

In this section we use the previous constructions to obtain a
particular tilting
bundle on any rational surfaces. 
We like 
to obtain a tilting bundle of small rank and a tilting bundle with a
quasi-hereditary endomorphism algebra. Note that we are not optimal
with our construction (compare for example with \cite{HP2}, Theorem
5.8), however, to avoid to many technical details and case by case
considerations we present one construction that 
works for any rational surface $X$.

It is convenient to start with 
a strongly exceptional sequence on a Hirzebruch surface $\varepsilon =
(\cO, \cO(P), \cO(Q + aP), \cO(Q + (a+1)P))$, where we can assume $a$
is sufficiently large. 

Then we chose a
sequence of blow ups, where we allow to blow up finitely many
different points in each step (so we use a slightly different notation
than in Section \ref{sectlbds})
$$
X= X^m \stackrel{\pi_{m}}{\lra} X^{m-1} \stackrel{\pi_{m-1}}{\lra}
\ldots X^2 \stackrel{\pi_{2}}{\lra} X^1 \stackrel{\pi_{1}}{\lra} X^0 =
X_0 = \Ff_a 
$$
of $X$ to a Hirzebruch surface. Note that we still have choices with
this notation. To make the morphisms unique (for a chosen surface $X$
with a fixed morphism $\pi$ to $\Ff_a$) we blow up in the first step
as many points as possible and proceed in this way. Thus if $x_v \in
X^i$ is a point not blown up under $\pi_{i+1}:X^{i+1} \lra X^i$ then
also its preimage in any $X^j$ for $j > i$ is not blown up. Moreover,
if $x_v \in X^i$ is blown up, we call $l(v) = i+1$ its level and denote
by $E_v$ the corresponding exceptional divisor in $X^{i+1}$.

Next we define the blow up graph $G$ as follows. Its vertices $G_0$
are the points $x_v$ in $X_i$ that are blown up under $\pi_i: X_{i+1}
\lra X_i$. There is an edge between $x_v$ and $x_w$, whenever $x_w \in
E_v$ (or vice versa). This way, we get a levelled graph, that is for
each edge $v$ --- $w$ we have $ |l(v)-l(w)| = 1$. The blow up graph is
precisely the Hasse diagram (Hasse graph) for the partial order
$\succ$ defined in Section \ref{sectlbds}. To construct an
exceptional sequence of line bundles on $X$ we also need the divsors
$R_v$ defined as the pull back of $E_v$ in $X$. Note that the strict
transform $\overline E_v$ of $E_v$ is an irreducible component of the
divisor 
$R_v$. For the 
self-intersection numbers on $X^{l(v)}$ we get $R_v^2 = -1 = E_v^2$
and on $X$ we obtain $R_v^2 = -1$ and $\overline
E_v^2 = -1 - a_v$ where $a_v$ is the number of points in $E_v$ blown
up under $\pi_{i+2}$.

To start with the construction and to avoid to many notation we assume
first $X$ is the recursiv blow up  of one point, so $\pi_i: X^i \lra
X^{i-1} $ is the blow up of one point $x_{i-1}$ on the exceptional
divisor $E_{i-1}$ for $i=1,\ldots,t = m$. With $R_i$ we denote the pull
back of $E_i$ to $X$. Then we consider the full exceptional sequence
$$
\varepsilon = (\cO, \cO(R_t),\cO(R_{t-1}),\ldots,\cO(R_2),
\cO(R_1),\cO(P),\cO(aP + Q), \cO((a+1)P + Q))
$$
that is obtained by recursive standard augmentations in the first
place. Using Lemma \ref{LExtOR} we obtain non-vanishing extension groups
$\mExt^1(\cO(R_i), \cO(R_j)) = \ck$ for all $i > j$. All other
extension groups vanish. Then we can define recursively vector bundles
$V_i$ via $V_1 = \cO(R_1)$ and
$$
0 \lra \cO(R_i) \lra V_i \lra V_{i-1} \lra 0
$$
with $\mExt^1(V_{i-1}, \cO(R_i)) = \ck$.
In this way we define vector bundles $V_i$ with non-trivial
endomorphism ring. In fact the direct sum of all $V_i$ has an
endomorphism ring isomorphic to the Auslander algebra of
$\ck[\alpha]/\alpha ^t$, a quasi-hereditary algebra  considered in
\cite{DR}, Section 7. 

\begin{Lemma}\label{Lexample1}
The coherent sheaf $V_i$ is a vector bundle of rank $i$ and it is
indecomposable with endomorphism ring isomorphic to
$\ck[\alpha]/\alpha ^i$. The direct sum $\oplus \cO \oplus \bigoplus_i V_i
\oplus \cO(P) \oplus \cO(aP + Q) \oplus
\cO((a+1)P + Q)$ is isomorphic to the
universal extension $\overline E$ of the exceptional sequence $\varepsilon$.
\end{Lemma}

{\sc Proof. }
Clearly $V_i$ is a vector bundle of rank $i$. To identify it with the universal
extension of $\cO(R_i)$ in the sequence $\varepsilon$ we need to
compute $\mExt^1_X(\cO(R_i),V_{i-1}) = \mExt^1_{X^i}(\cO(E_i),V_{i-1}) =
k$. This can be shown recursively, just apply $\mHom(\cO(E_i),-)$ to the
defining sequence and we obtain $\mHom(\cO(E_i), \cO(R_j)) = k =
\mHom(\cO(E_i),V_j)$ and $\mExt^1(\cO(E_i), \cO(R_j)) = k =
\mExt^1(\cO(E_i),V_j) $ for all $j < i$. This follows directly from
the equivalence in Theorem \ref{Tgeneralequivalence} for the
exceptional sequence  
$$
(\cO(R_t),\cO(R_{t-1}),\ldots,\cO(R_2), \cO(R_1)).
$$
Note that the quasi-hereditary algebra $A$ for this
sequence is the Auslander algebra of $\ck[\alpha]/\alpha ^i$
considered in \cite{DR}, Section 7. However, one might check the claim
also directly using the defining exact sequences and the vanishing of $\mExt^2$.
\hfill $\Box$

Now we consider the general case, let $X$ be any rational surface
together with a sequence of blow ups from a Hirzebruch surface. Then
the $\mExt$--blocks correspond to the points in $X^0 = \Ff_a$ that are
blown up under $\pi_1$. Consequently, the $\mExt$--blocks correspond
to the connected components in the blow up graph, together with the
four bundles $\cO, \cO(P), \cO(aP + Q)$ and $\cO((a+1)P + Q)$ we
started with. To be precise, for our exceptional sequence the
$\mExt$--graph of $(\cO(R_t),\cO(R_{t-1}),\ldots,\cO(R_2), \cO(R_1))$
contains the blow up graph, by Lemma \ref{LExtOR}. Moreover, they
have the same connected components. Note that we get non-trivial
universal extensions only between two objects in the same connected
component. For any point $x_j \in X^{l(j)-1}$ blown up 
under some morphism $\pi_{l(j)}$ we define the universal extension $V_j$
of $\cO(R_j)$ by all bundles $\cO(R_i)$ with $x_j$ is blown down to
$x_i$ under some composition of the maps $\pi$. This bundle $V_i$,
according to Lemma \ref{Lexample1}, is a direct summand of the universal extension
$\overline E$ for the exceptional sequence $\varepsilon$. The
arguments above for the particular case also apply here, so we get a
tilting bundle on $X$ satisfying the following properties.

\begin{Thm}
Let $\overline E$ be the universal extension of the full exceptional
sequence $\varepsilon$ above, then the direct summmands of $\overline
E$ are isomorphic to the vector bundles $V_i$. In particular, $
\cO \oplus \bigoplus_i V_i \oplus \cO(P) \oplus \cO(aP + Q) \oplus
\cO((a+1)P + Q)$ is a tilting bundle on $X$. If $x_i \in X^j$ then $\rk
V_i = j$. Thus $\overline E$ consists of vector bundles of rank at
most $t$. Moreover, $\mHom(V_i,V_j) \not= 0$ precisely when $x_i$ blows down to
$x_j \in X^l$. In this case $\mHom(V_i,V_j)$ is $l$--dimensional.
\end{Thm}

{\sc Proof. }
First note that the universal extension of $\cO(R_i)$ with respect to
$\varepsilon$ coincides with the universal extension for the
exceptional sequence consisting of all $\cO(R_j)$ with $x_i$ is mapped
to $xj$ under the composition $X^{l(i)} \lra X^{l(j)}$, since all
other extension groups vanish. Thus we can apply the lemma above to
show that $V_i$ is indecomposable of rank $l(i)$ and $\overline E$
consists of the four line bundles on $\Ff_a$ and the bundles
$V_i$. Consequently, $\overline
E$ are isomorphic to the vectorbundles $V_i$. In particular, $
\cO \oplus \bigoplus_i V_i \oplus \cO(P) \oplus \cO(aP + Q) \oplus
\cO((a+1)P + Q)$ is a tilting bundle on $X$ with $\rk V_i$ equals the
level of $x_i$. Finally, the claim for the Hom--groups can be shown by
induction using the defining exact sequences.
\hfill $\Box$

{\sc Remark.}
We have chosen a simple example to construct at least one particular tilting
bundle. In fact we have many other choices. First, we can use
different projections to different Hirzebruch surfaces. Then we can choose
the position of any standard augmentation and, eventually, we can
chose to perform either extensions or coextensions. Moreover, we do
not need to start with line bundles, in fact also the structure sheaf
on any $(-1)$-curve can be used, since it is exceptional. However,
apart from this a construction of other exceptional sheaves becomes
more technical and the computation of the extensions groups might be
more difficult as well. Thus line bundles are a nice, but not the only choice.

\medskip

{{\small Lutz Hille}\\
{\small Mathematisches Institut}\\
{\small Universit\"at M\"unster}\\
{\small Einsteinstr.~62}\\
{\small D--48149 M\"unster}\\
{\small Germany}\\
{\small E-mail: lutz.hille@uni-muenster.de}\\
{\small 
}}
\medskip

{{\small Markus Perling}\\
{\small Ruhr-Universit\"at Bochum }\\
{\small Fakult\"at f\"ur Mathematik }\\
{\small Universit\"atsstrasse 150 }\\
{\small D--44780 Bochum }\\
{\small Germany}\\
{\small E-mail: markus.perling@rub.de}\\
{\small
}}

\end{document}